# Discrete Groups, Grothendieck Rings and Families of Finite Subgroups

ALEJANDRO ADEM

Dedicated to Samuel Gitler on occasion of his sixtieth birthday.

## 0. Introduction

Understanding the representation theory of an infinite discrete group is usually very difficult. Nevertheless, there are certain classes of discrete groups (such as arithmetic groups) which satisfy algebraic properties similar to those of a finite group. For example, an arithmetic group has only a finite number of elements of finite order up to conjugacy, and the same is true for other large classes of groups such as the mapping class groups, the outer automorphism group of a free group, hyperbolic groups, etc.

In this paper we look at the representation theory which is determined by *finite subgroups*. More precisely, if $\Gamma$ is a discrete group with a finite number of finite subgroups up to conjugacy (the FCP condition) we introduce

DEFINITION *If $\mathcal{F}$ is a family of finite subgroups, and $R$ is a commutative ring, then*
$$G_{\mathcal{F}}(R[\Gamma]) = \varprojlim_{H \in \mathcal{F}} G(R[H]),$$
*where $G(R[H])$ denotes the usual Grothedieck ring of finitely generated $R[H]$–modules.*

After providing preliminaries in section 1, we prove some basic facts about the rings $G_{\mathcal{F}}(R[\Gamma])$ which follow easily from the corresponding results for finite groups. We summarize them below (notation as above):
(1) If $S$ is a commutative ring with 1 satisfying the minimal conditions on ideals then $G_{\mathcal{F}}(S[\Gamma])$ is free abelian of finite rank.

1991 *Mathematics Subject Classification.* 20C07, 55R40, 19A22.
Partially supported by an NSF grant, an NSF Young Investigator Award and the ETH–Zürich.
This paper is in final form and no version of it will be submitted for publication elsewhere.

(2) $rk_{\mathbb{Z}} G_{\mathcal{F}}(\mathbb{Q}[\Gamma]) = C_{\mathcal{F}}(\Gamma)$, the number of distinct conjugacy classes of cyclic subgroups in $\mathcal{F}$.

(3) $rk_{\mathbb{Z}} G_{\mathcal{F}}(\mathbb{C}[\Gamma]) = n_{\mathcal{F}}(\Gamma)$, the number of distinct conjugacy classes of $\gamma \in \Gamma$ such that $<\gamma> \in \mathcal{F}$.

(4) For any field $\mathbb{K}$ of characteristic zero, $G_{\mathcal{F}}(\mathbb{K}[\Gamma]) = \mathbb{Z}$ if and only if $\mathcal{F} = \{\{1\}\}$.

(5)
$$dim_{\mathbb{Q}} G_{\mathcal{F}}(\mathbb{Z}[\Gamma]) \otimes \mathbb{Q} \geq C_{\mathcal{F}}(\Gamma).$$

In §3 we compute $G_{\mathcal{F}}(\mathbb{Q}[\Gamma])$, $G_{\mathcal{F}}(\mathbb{C}[\Gamma])$ for several interesting examples including amalgamated products, $SL_3(\mathbb{Z})$, $GL_{p-1}(\mathbb{Z})$ and a few of the Bianchi groups. These examples illustrate the computability of our construction.

Let $G^f(R[\Gamma])$ denote the Grothendieck ring constructed from $R[\Gamma]$–modules which are finitely generated over $R$, and define

$$\overline{G}_{\mathcal{F}}(R[\Gamma]) = G^f(R[\Gamma]) \bigg/ \bigcap_{H \in \mathcal{F}} \ker \mathrm{res}_H^{\Gamma}.$$

In section 2 we discuss how this ring compares to the previous one, measured by an inclusion

$$j_{\mathcal{F}} : \overline{G}_{\mathcal{F}}(R[\Gamma]) \to G_{\mathcal{F}}(R[\Gamma]).$$

If $\Gamma$ is a finite group, the usual induction–restriction arguments can be used to show that the map above is a rational equivalence. In contrast, this is *not true* for an arbitrary discrete group; the two rings may have different rank. This important discrepancy indicates that our construction will be of *independent interest*, irrespective of what the actual Grothendieck group associated to $R[\Gamma]$ looks like.

However this intermediate ring is useful for some specific questions relating to representations of $\Gamma$ restricted to their finite subgroups. To illustrate this we introduce a cohomological index

$$\overline{G}_{\mathcal{F}}(K[\Gamma]) \longrightarrow \mathbb{Z}$$

for certain groups $\Gamma$, showing that this ring is a natural object for questions on Euler characteristics and related invariants.

One of the central problems in algebraic $K$–theory has been to compute the Grothendieck group of projective modules $K_0(\mathbb{Z}[\Gamma])$ (see [F-J]). In particular, if $\Gamma$ is torsion–free, it is conjectured that $K_0(\mathbb{Z}[\Gamma]) \cong \mathbb{Z}$, and this has been verified in many instances using topological methods. In §4 we apply our methods to analyze $K_0(\mathbb{K}[\Gamma])$ when $\Gamma$ has torsion. Using a construction introduced by Loday [L] we prove a converse to this, namely

PROPOSITION. *If $\Gamma$ is FCP and $\mathbb{K}$ is a field of characteristic zero, then for $\mathcal{F}$, the family of all finite subgroups in $\Gamma$, we have*

$$dim_{\mathbb{Q}} K_0(\mathbb{K}[\Gamma]) \otimes \mathbb{Q} \geq dim_{\mathbb{Q}} (im \, j_{\mathcal{F}}) \otimes \mathbb{Q}.$$

In particular, if $\Gamma$ is virtually torsion–free (see section 1 for definition), we obtain the converse of the conjecture mentioned above over the complex numbers: if $K_0(\mathbb{C}[\Gamma]) \cong \mathbb{Z}$, then $\Gamma$ is torsion–free.

In the final section we show how our construction is useful in equivariant topological K–theory. If $G$ is a finite group and $Y$ is a finite $G$–CW complex, we construct a limit $\mathcal{R}_G(Y)$ using the representation rings of the isotropy subgroups and fixed–point data (see section 5 for precise definitions) as well as a map induced by the restrictions $\varphi : K_G^0(Y) \longrightarrow \mathcal{R}_G(Y)$ for which we prove

PROPOSITION. *If $Y$ is a finite $G$–CW complex, then the ring map $\varphi$ has nilpotent kernel and is a rational surjection.*

As a special case of this, under appropriate hypotheses we obtain $G_{\mathcal{F}}(\mathbb{C}[\Gamma]) \otimes \mathbb{Q}$ as a quotient of the equivariant K–theory of an admissible $\Gamma$–complex. In fact this paper was motivated by this topological application, which appeared in [A].

The simple algebraic constructions we use seem to be of some interest in their own right, especially given the fact that they record an important difference between finite groups and examples such as arithmetic groups. Moreover the applications above indicate that the techniques described in this paper may be of some use in topological K–theory, algebraic $K$–theory, group cohomology, and character theory. Our main goal in writing this paper has been to provide an accessible account of this, leaving many questions which arise in this context to the interested reader.

It is a pleasure to dedicate this paper to Samuel Gitler, to whom I owe a debt of gratitude for both his mathematical inspiration and friendship.

## 1. Preliminaries

In this section we will briefly discuss the class of groups to which our methods will apply. Let $\Gamma$ denote a discrete group.

DEFINITION 1.1. *A family $\mathcal{F}$ of subgroups of $\Gamma$ is a non–empty collection of subgroups which is closed under conjugation and taking subgroups.*

Note in particular that $\{1\}$ will belong to any family. By definition, $\Gamma$ will act on a family $\mathcal{F}$ via conjugation.

DEFINITION 1.2. *Given a discrete group $\Gamma$ we say that it satisfies the finite conjugacy property (FCP) if there are only finitely many finite subgroups in $\Gamma$ up to conjugacy.*

If $\Gamma$ satisfies FCP, and $\mathcal{F}$ is any family of finite subgroups in $\Gamma$, then we denote by $n_{\mathcal{F}}(\Gamma)$ the number of distinct conjugacy classes of $\gamma \in \Gamma$ such that $<\gamma> \in \mathcal{F}$. This will be a finite integer.

¿From now on we will only be dealing with families of *finite* subgroups, and with groups which satisfy the FCP condition.

EXAMPLE 1.1. We now list certain specific families we will be interested in:

(i)     $\mathcal{F}(\Gamma)$    =    family of all finite subgroups in $\Gamma$,

(ii)     $\mathcal{F}_p(\Gamma)$    =    family of all finite $p$– subgroups in $\Gamma$,

(iii)     $\mathcal{F}_{p'}(\Gamma)$    =    family of all finite $p'$– subgroups in $\Gamma$,

(iv)     $\mathcal{C}(\Gamma)$    =    family of all cyclic finite subgroups in $\Gamma$.

To minimize our notation, we will denote $n_{\mathcal{F}(\Gamma)}(\Gamma)$ by $n(\Gamma)$, $n_{\mathcal{F}_p(\Gamma)}(\Gamma)$ by $n_p(\Gamma)$, and $n_{\mathcal{F}_{p'}(\Gamma)}(\Gamma)$ by $n_{p'}(\Gamma)$. Note that if $\mathcal{F}$ is any family, we can define a new family from it,

$$\mathcal{C}(\mathcal{F}) = \text{family of all cyclic subgroups in } \mathcal{F}$$

and that clearly $n_{\mathcal{F}}(\Gamma) = n_{\mathcal{C}(\mathcal{F})}(\Gamma)$.

Examples of groups satisfying FCP include finite groups, arithmetic groups, mapping class groups, and many other interesting classes (see [A]).

For later purposes, we will need an additional hypothesis on $\Gamma$. We recall

DEFINITION 1.3. *If $\Gamma$ is a discrete group, it is said to be virtually torsion–free (VTF) if it contains a torsion–free subgroup $\Gamma'$ of finite index in $\Gamma$.*

The examples mentioned previously are all virtually torsion–free. For example if $\Gamma = SL_n(\mathbb{Z})$, we can take $\Gamma' = \Gamma(p)$, the level $p$ congruence subgroup for $p$ an odd prime. Note that we can always choose $\Gamma'$ to be normal in $\Gamma$.

Next, we recall the concept of the Grothendieck group associated to a ring. Let $A$ be a ring; the Grothendieck group $G(A)$ of finitely generated left $A$–modules is defined (see [Sw1]) to be the group with generators $[M]$, $M$ an $A$–module, and relations

$$[M] = [M'] + [M'']$$

whenever there is an exact sequence

$$0 \longrightarrow M' \longrightarrow M \longrightarrow M'' \longrightarrow 0$$

of $A$–modules. Addition is induced by direct sum.

If $R$ is a commutative ring, and $A$ is an $R$–algebra which is finitely generated and projective as an $R$–module, then we may do the same construction as above but using only modules which are finitely generated and projective over $R$; we denote this by $G^R(A)$. If we use only projective modules, we obtain $K_0(A)$. There is a natural map

$$G^R(A) \xrightarrow{P_A} G(A)$$

which is an isomorphism if $R$ is a regular commutative ring. The particular case we will consider is $A = R[\pi]$, where $\pi$ is a discrete group, and $R$ is a regular commutative ring (we will make this assumption from now on). By the above, $G^R(R[\pi])$ is an associative, commutative ring with unit, and if $\pi$ is finite, then $P_A$ is an isomorphism. In addition, a homomorphism $R \longrightarrow R'$ of regular commutative rings induces a map $G(R[\pi]) \longrightarrow G(R'[\pi])$.

## 2. Limits of Grothendieck Groups

We assume now that $\Gamma$ is a discrete group satisfying the FCP condition and that $R$ is a regular, commutative ring.

DEFINITION 2.1. *The $\mathcal{F}$–Groethendieck ring of $R[\Gamma]$ is*

$$G_\mathcal{F}(R[\Gamma]) = \lim_{H \in \mathcal{F}} G(R[H])$$

The limit on the right can be described explicitly as follows. Conjugation by $\gamma \in \Gamma$ induces a group homomorphism $C_\gamma : H \longrightarrow \gamma H \gamma^{-1}$, and hence, a ring homomorphism

$$C_\gamma^* : G(R[\gamma H \gamma^{-1}]) \longrightarrow G(R[H]).$$

Similarly, the inclusion $H \hookrightarrow K$ induces a restriction map

$$\mathrm{res}_H^K : G(R[K]) \longrightarrow G(R[H]).$$

Then the limit in 2.1 is, by definition, the subset of sequences

$$(x_H) \in \prod_{H \in \mathcal{F}} G(R[H])$$

such that
  (1) if $H \subseteq K$, then $\mathrm{res}_H^K(x_K) = x_H$
  (2) if $K = \gamma H \gamma^{-1}$, then $C_\gamma^*(x_K) = x_H$.

As the families we consider have a finite number of elements up to conjugacy, we can in fact identify this limit with a subset of the direct sum

$$\bigoplus_{\substack{(H) \\ H \in \mathcal{F}}} G(R[H])^{w_\Gamma(H)}$$

satisfying the stability condition (1), where $w_\Gamma(H) = N_\Gamma(H)/C_\Gamma(H)$ is a finite group, and $(H)$ denotes the conjugacy class represented by $H$.

Given our hypotheses, $G_\mathcal{F}(R[\Gamma])$ is a unitary, associative and commutative ring of *finite type*. We will examine some of its general properties before computing specific examples. Most of them follow easily from the corresponding facts for finite groups (see [Sw1], [Sw2]).

To begin we have

PROPOSITION 2.1. *Let $S$ be any commutative ring with unit satisfying the minimal condition on ideals. Then, for any family $\mathcal{F}$ of finite subgroups in $\Gamma$, $G_\mathcal{F}(S[\Gamma])$ is a free abelian group of finite rank.*

Next, note that if $R \longrightarrow R'$ is a map of commutative rings, then, by the considerations in §1, there is an induced map

$$G_{\mathcal{F}}(R[\Gamma]) \longrightarrow G_{\mathcal{F}}(R'[\Gamma]).$$

In particular, if $\mathbb{K}$ is the quotient field of $R$ (an integral domain), we have a ring homomorphism

$$G_{\mathcal{F}}(R[\Gamma]) \xrightarrow{\psi_{\mathcal{F}}} G_{\mathcal{F}}(\mathbb{K}[\Gamma]).$$

Recall that if $H$ is a finite group, and $\mathbb{K}$ is a field of characteristic zero, the character ring of $\mathbb{K}[H]$ is the set of all integral linear combinations of characters of representations of $H$ over $\mathbb{K}$, denoted by $C(\mathbb{K}[H])$. The character map

$$\begin{aligned} G(\mathbb{K}[H]) &\longrightarrow C(\mathbb{K}[H]) \\ A &\longmapsto \chi_A \end{aligned}$$

induces a ring isomorphism, and hence, (as these maps are clearly compatible) we have

PROPOSITION 2.2. *The character maps induce an isomorphism*

$$G_{\mathcal{F}}(\mathbb{K}[\Gamma]) \cong \varprojlim_{H \in \mathcal{F}} C(\mathbb{K}[H])$$

*for any family $\mathcal{F}$ of finite subgroups in $\Gamma$.*

Hence $G_{\mathcal{F}}(\mathbb{K}[\Gamma])$ is clearly of finite rank, and in principle can be computed. In analogy with classical representation theory, we would like to compute the rank of $G_{\mathcal{F}}(\mathbb{K}[\Gamma])$ in terms of subgroup data in $\Gamma$ related to $\mathcal{F}$. Denote by $C_{\mathcal{F}}(\Gamma)$ the number of distinct conjugacy classes of cyclic subgroups in $\mathcal{F}$; we have

PROPOSITION 2.3. $G_{\mathcal{F}}(\mathbb{Q}[\Gamma])$ *is free abelian of rank* $C_{\mathcal{F}}(\Gamma)$.

PROOF. For a finite group $H$, let $\mathcal{C}(H)$ denote the family of its cyclic subgroups. Then $G(\mathbb{Q}[H]) \otimes \mathbb{Q}$ can be identified with class functions

$$\mathcal{C}(H) \longrightarrow \mathbb{Q}.$$

Now we can also identify $\varinjlim_{H \in \mathcal{F}} \mathcal{C}(H)$ with the family $\mathcal{C}(\mathcal{F})$. Hence, there is an isomorphism

$$\mathrm{Map}_{\Gamma}(\mathcal{C}(\mathcal{F}), \mathbb{Q}) \cong \varprojlim_{H \in \mathcal{F}} \mathrm{Map}_H(\mathcal{C}(H), \mathbb{Q}).$$

The term on the right is isomorphic to $G_{\mathcal{F}}(\mathbb{Q}[\Gamma]) \otimes \mathbb{Q}$, hence showing that the original ring is free abelian of rank

$$dim_{\mathbb{Q}} \mathrm{Map}_{\Gamma}(\mathcal{C}(\mathcal{F}), \mathbb{Q}) = C_{\mathcal{F}}(\Gamma).$$

□

In a similar way one can show

PROPOSITION 2.4. *$G_{\mathcal{F}}(\mathbb{C}[\Gamma])$ is free abelian of rank $n_{\mathcal{F}}(\Gamma)$.*

Note in particular that if $\Gamma$ is finite and $\mathcal{F} = \mathcal{F}_{p'}(\Gamma)$, then $G_{\mathcal{F}}(\mathbb{C}[\Gamma])$ is a ring which is free abelian of rank $n_{p'}(\Gamma)$, just like the usual ring of modular characters. Upon tensoring with a large enough field they clearly coincide, but it is not clear if they are equal. If we take $\mathcal{F} = \mathcal{F}_p(\Gamma)$ then we obtain a ring of rank the number of conjugacy classes of elements of order a power of $p$ in $\Gamma$. Such a ring would appear to be related to $p$–adic representations of the group.

In fact, following the usual methods in character theory ([Sw2]), one can prove similar results for any field of characteristic zero. We obtain a general result

PROPOSITION 2.5. *If $\mathbb{K}$ is a field of characteristic zero, and $\mathcal{F}$ is any family of finite subgroups in $\Gamma$, then $G_{\mathcal{F}}(\mathbb{K}[\Gamma]) = \mathbb{Z}$ if and only if $\mathcal{F} = \{\{1\}\}$.*

Recall that if $\mathcal{F}$ is a family, $C(\mathcal{F})$ denotes the subfamily of all the cyclic subgroups in $\mathcal{F}$. We have a natural map

$$G_{\mathcal{F}}(R[\Gamma]) \to G_{\mathcal{C}(\mathcal{F})}(R[\Gamma]).$$

Let $R = \mathbb{K}$, a field of characteristic zero. Then, as characters are determined on cyclic subgroups, this map is injective, and in fact we obtain

PROPOSITION 2.6. *For any field $\mathbb{K}$ of characteristic zero,*

$$G_{\mathcal{F}}(\mathbb{K}[\Gamma]) \otimes \mathbb{Q} \cong G_{\mathcal{C}(\mathcal{F})}(\mathbb{K}[\Gamma]) \otimes \mathbb{Q}. \qquad \Box$$

Families of cyclic subgroups are particularly nice for computing $G_{\mathcal{F}}(\mathbb{Z}[\Gamma])$. Using the fact that for any family $\mathcal{C}$ of cyclic subgroups, and $H \in \mathcal{C}$ the maps $G(\mathbb{Z}[H]) \longrightarrow G(\mathbb{Q}[H])$ split in a compatible way [Sw3], we obtain

PROPOSITION 2.7. *Let $\mathcal{C}$ be a family of cyclic subgroups of $\Gamma$. Then the natural map*

$$G_{\mathcal{C}}(\mathbb{Z}[\Gamma]) \longrightarrow G_{\mathcal{C}}(\mathbb{Q}[\Gamma])$$

*is a split surjection of rings.*

Combining these results we obtain

PROPOSITION 2.8. *For any family $\mathcal{F}$ of finite subgroups,*

$$dim_{\mathbb{Q}} \, G_{\mathcal{F}}(\mathbb{Z}[\Gamma]) \otimes \mathbb{Q} \geq C_{\mathcal{F}}(\Gamma). \qquad \Box$$

More generally it is known [Sw3] that for a finite group $H$, if $\mathbb{K}$ is the field of fractions of $R$, a regular, Noetherian integral domain, then the map $G(R[H]) \to G(\mathbb{K}[H])$ is a surjection, but not necessarily split as a map of rings. Hence there is no clear reason to expect the map

$$G_{\mathcal{F}}(R[\Gamma]) \to G_{\mathcal{F}}(\mathbb{K}[\Gamma])$$

to be surjective for an arbitrary family $\mathcal{F}$ of finite subgroups.

The basic results in this section indicate that for any discrete group $\Gamma$ satisfying the FCP condition our construction $G_{\mathcal{F}}(R[\Gamma])$ is an accessible algebraic invariant which can be used to concentrate on representation–theoretic problems both away from $p$ and at $p$. We will see later on that it is computable in many interesting cases.

Our previous construction is related to a particular Grothendieck group of $R[\Gamma]$–modules. We introduce

DEFINITION 2.2. *Let $R$ be a commutative ring and $\Gamma$ a discrete group, we define $G^f(R[\Gamma])$ to be the Grothendieck group of $R[\Gamma]$–modules which are finitely generated over $R$.*

Given that $\Gamma$ is FCP, the restrictions induce a map $G^f(R[\Gamma]) \longrightarrow \prod G(R[H])$ which, using the natural inclusion of stable elements in the product, can be factorized as

$$G^f(R[\Gamma]) \longrightarrow G_{\mathcal{F}}(R[\Gamma])$$
$$\searrow \qquad \nearrow$$
$$\prod_{H\in\mathcal{F}} G(R[H])$$

DEFINITION 2.3.
$$\overline{G}_{\mathcal{F}}(R[\Gamma]) = \mathrm{im}(\prod_{H\in\mathcal{F}} \mathrm{res}) \qquad \square$$

Note that we have a natural inclusion $j_{\mathcal{F}} : \overline{G}_{\mathcal{F}}(R[\Gamma]) \to G_{\mathcal{F}}(R[\Gamma])$; the key difference between these rings is that one of them is a *quotient* of a representation ring,

$$\overline{G}_{\mathcal{F}}(R[\Gamma]) \cong G^f(R[\Gamma]) \Big/ \bigcap_{H\in\mathcal{F}} \ker \mathrm{res}_H^{\Gamma}.$$

If $\Gamma$ is finite, then standard induction–restriction arguments can be used to show that $j_{\mathcal{F}}$ becomes a surjection after tensoring with the rationals. However, it has been pointed out to the author [Se2] that this does not hold for arbitrary discrete groups. This is an important point, as it indicates that our construction distinguishes conjugacy classes of elements of finite order, whereas the actual representation ring *does not*. This is a significant departure from the behaviour for finite groups, and would seem to indicate that the rings $G_{\mathcal{F}}(R[\Gamma])$ are of interest *independently* of actual representations. On the other hand the ring $\overline{G}_{\mathcal{F}}(R[\Gamma])$ can be useful in dealing with algebraic and topological questions determined by actual representations of $\Gamma$ restricted to finite subgroups.

Let us assume that $\Gamma$ is VTF, and that $\Gamma' \subseteq \Gamma$ is a normal subgroup of finite index in $\Gamma$ with quotient $G = \Gamma/\Gamma'$ (a finite group), and $\pi : \Gamma \twoheadrightarrow G$ the projection map. We have an induced map of rings $G(R[G]) \to \overline{G}_{\mathcal{F}}(R[\Gamma])$, the image of which we will denote by $\{\pi_*(G(R[G]))\}$. This can sometimes be used to give a lower bound on the rank of the ring $\overline{G}_{\mathcal{F}}(R[\Gamma])$. For example, if $R = \mathbb{C}$, let $n_G$ be equal to the number of distinct $G$–conjugacy classes of elements $g \in G$ of the form $g = \pi(\gamma)$, where $\gamma \in \Gamma$ is an element of finite order. Then, using character theory arguments, one can prove that in fact the rank of $\{\pi_*(G(\mathbb{C}[G])\}$ is precisely $n_G$, and hence we obtain (for $\mathcal{F}$ the family of all finite subgroups in $\Gamma$)

PROPOSITION 2.9.

$$dim_{\mathbb{Q}} \overline{G}_{\mathcal{F}}(\mathbb{C}[\Gamma]) \otimes \mathbb{Q} \geq max\ n_G$$

*as $G$ ranges over all finite quotients of $\Gamma$ with torsion–free kernel.*

To illustrate the usefulness of the intermediate ring described above, we now show how it is the natural domain for a cohomological index closely related to the Euler characteristic of $\Gamma$ (see [Br] for background on this). In the remainder of this section we will make the additional assumption that there exists a finite dimensional, acyclic $\mathbb{Z}[\Gamma]$–complex $C_*$ such that for each $i \geq 0$,

$$C_i \cong \bigoplus_{k=1}^{m_i} \mathbb{Z} \bigotimes_{H_{i,k}} \mathbb{Z}[\Gamma],$$

where each $H_{i,k}$ is a *finite* subgroup of $\Gamma$. For geometric reasons, this assumption holds for arithmetic groups, amalgamated products of finite groups, and many other interesting examples [A]. In this section we denote $\mathcal{F} = \mathcal{F}(\Gamma)$, the family of all finite subgroups in $\Gamma$. Let $m(\Gamma) = $ l.c.m. $\{|H| \mid H \in \mathcal{F}\}$. We introduce

DEFINITION 2.4. *Let $K$ be a field of characteristic prime to $m(\Gamma)$, and $M$ a $K[\Gamma]$–module of finite rank over $K$. The index of $M$ is, by definition, the integer*

$$I_K(M) = \sum_{i=0}^{\infty} (-1)^i \dim_K H^i(\Gamma, M) \qquad \square$$

PROPOSITION 2.10. *With the previous hypotheses, we have*

*(1) $I_K(M)$ is a finite integer.*

*(2) If $M_1\big|_S \cong M_2\big|_S$ for every $S \in \mathcal{F}$, then $I_K(M_1) = I_K(M_2)$.*

*(3) If $0 \longrightarrow M'' \longrightarrow M \longrightarrow M' \longrightarrow 0$ is a short exact sequence, then*

$$I_K(M) = I_K(M') + I_K(M'').$$

PROOF. Let $N = \dim\ C_*$ (as a chain complex); then we have

$$I_K(M) = \sum_{i=0}^{N} (-1)^i \dim_K \left( \bigoplus_{k=1}^{m_i} M^{H_{i,k}} \right).$$

Hence, (1), and (2) follow. As for (3), this is a conseqence of the fact that taking invariants under any $Q \in \mathcal{F}$ is exact, given our hypothesis on char$(K)$. $\square$

By Definition 2.2, we see that the index induces a well defined additive map

$$I_K : \overline{G}_\mathcal{F}(K[\Gamma]) \longrightarrow \mathbb{Z}.$$

Note that $I_K(K) = \chi(B\Gamma)$, the topological Euler characteristic of $\Gamma$. Assume that $\Gamma$ is virtually torsion–free, and let $\Gamma' \subseteq \Gamma$ be a normal, torsion–free subgroup with finite quotient $G = \Gamma/\Gamma'$. Let $K[G]$ denote the group algebra of $G$ over $K$, and suppose it decomposes into a sum of indecomposable modules,

$$K[G] \cong \bigoplus_{r=1}^{\ell(G)} S_r.$$

We can think of this as an isomorphism of $K[\Gamma]$–modules using the projection $\Gamma \twoheadrightarrow G$, and hence, we obtain

$$I_K(K[G]) = \sum_{r=1}^{\ell(G)} I_K(S_r).$$

On the other hand, $H^*(\Gamma, K[G]) \cong H^*(\Gamma', K)$, and so

$$I_K(K[G]) = \chi(B\Gamma').$$

Combining these facts, we obtain that $H^*(\Gamma', K) \cong \bigoplus H^*(\Gamma, S_r)$, and that

$$\chi(B\Gamma') = \sum_{r=1}^{\ell(G)} I_K(S_r).$$

Note that if $\Gamma$ is torsion–free, then $I_K(M) = \dim_K M \cdot \chi(B\Gamma)$ for all $M$, and hence $\operatorname{im} I_K \cong \chi(B\Gamma) \cdot \mathbb{Z} \subseteq \mathbb{Z}$ in this case. In the general case it is not clear what $\operatorname{im} I_K$ must be.

As an application of this invariant, we use it to give a proof of a result due to Brown [Br].

PROPOSITION 2.11. *If $\Gamma' \subseteq \Gamma$ is a torsion–free subgroup of finite index, then $(m(\Gamma)/[\Gamma : \Gamma']) \cdot \chi(B\Gamma')$ is an integer.*

PROOF. We already pointed out that $\chi(B\Gamma') = I_K(K \otimes_{\Gamma'} K[\Gamma])$. Now, for any $H \in \mathcal{F}$,

$$K \bigotimes_{\Gamma'} K[\Gamma] \bigg|_H \cong \frac{[\Gamma : \Gamma']}{|H|} K[H],$$

and so

$$\operatorname{res}_H^\Gamma (m(\Gamma) \cdot K \bigotimes_{\Gamma'} K[\Gamma]) = [\Gamma : \Gamma'] \, (q_H \cdot K[H])$$

where $q_H = m(\Gamma)/|H|$. We can compute explicitly from this, to obtain

$$m(\Gamma) \cdot I_K(K \bigotimes_{\Gamma'} K[\Gamma]) = \left( \sum_{i=0}^{N} (-1)^i \left( \sum_{k=1}^{m_i} q_{H_{i,k}} \right) \right) [\Gamma : \Gamma'],$$

and so
$$m(\Gamma) \Big/ [\Gamma : \Gamma'] \cdot \chi(B\Gamma') \in \mathbb{Z}.$$

□

## 3. Calculations over the Rational and Complex Numbers

In this section we will calculate $G_{\mathcal{F}}(\mathbb{Q}[\Gamma])$, $G_{\mathcal{F}}(\mathbb{C}[\Gamma])$ for several interesting examples, illustrating the computability of these invariants.

We start by considering the amalgamated product of two finite groups along a common subgroup, $\Gamma = G_1 *_H G_2$. Then we have an exact sequence

$$(3.1) \quad 0 \longrightarrow G_{\mathcal{F}}(R[\Gamma]) \longrightarrow G(R[G_1]) \oplus G(R[G_2]) \xrightarrow{\operatorname{res}_H^{G_1} - \operatorname{res}_H^{G_2}} G(R[H])$$

for $\mathcal{F}$ = family of all finite subgroups in $\Gamma$. We compute some specific examples for $R = \mathbb{Q}$.

EXAMPLE 3.1. (1) $\Gamma = \mathbb{Z}/p * \mathbb{Z}/q$  $p, q$  distinct primes
In this case (3.1) becomes

$$0 \longrightarrow G_{\mathcal{F}}(\mathbb{Q}[\Gamma]) \longrightarrow G(\mathbb{Q}[\mathbb{Z}/p]) \oplus G(\mathbb{Q}[\mathbb{Z}/q]) \longrightarrow \mathbb{Z} \longrightarrow 0.$$

Recall, that, for $p$ a prime,

$$G(\mathbb{Q}[\mathbb{Z}/p]) \cong \mathbb{Z}[w] \Big/ w^2 - (p-2)w - (p-1).$$

The ring $G_{\mathcal{F}}(\mathbb{Q}[\Gamma])$ can be computed from this information. We do the particular case when $p = 3$ and $q = 5$:

$$G(\mathbb{Q}[\mathbb{Z}/3]) \cong \mathbb{Z}[x] \Big/ x^2 - x - 2,, \quad G(\mathbb{Q}[\mathbb{Z}/5]) \cong \mathbb{Z}[x] \Big/ y^2 - 3y - 4,$$

and then,

$$G(\mathbb{Q}[\Gamma]) \cong \mathbb{Z}[z] \Big/ z^3 - 4z^2 - z + 1.$$

Note that the rank is three in this case, corresponding to the fact that there are three cyclic subgroups of finite order in $\Gamma$, up to conjugacy.

(2) $\Gamma = SL_2(\mathbb{Z})$

This group can be expressed as the amalgamated product $\mathbb{Z}/4 *_{\mathbb{Z}/2} \mathbb{Z}/6$, hence we have the short exact sequence

$$0 \longrightarrow G_{\mathcal{F}}(\mathbb{Q}[\Gamma]) \longrightarrow G(\mathbb{Q}[\mathbb{Z}/4]) \oplus G(\mathbb{Q}[\mathbb{Z}/6]) \xrightarrow{\operatorname{res}_1 - \operatorname{res}_2} G(\mathbb{Q}[\mathbb{Z}/2]).$$

Now

$$G(\mathbb{Q}[\mathbb{Z}/2]) \cong \mathbb{Z}[z] \Big/ z^2 - 1,$$

$$G(\mathbb{Q}[\mathbb{Z}/4]) \cong \mathbb{Z}[x,y] \Big/ \begin{array}{l} x^2 - 1, xy - y, \\ y^2 - 2x - 2 \end{array}$$

and
$$G(\mathbb{Q}[\mathbb{Z}/6]) \cong \mathbb{Z}[s,t] \Big/ \begin{array}{l} s^2 - s - 1, \\ t^2 - 1. \end{array}$$

The restriction maps take the values

$$\operatorname{res}_1(x) = 1, \quad \operatorname{res}_1(y) = 2z, \quad \operatorname{res}_2(s) = 2, \quad \operatorname{res}_2(t) = z,$$

from which we deduce

$$G_{\mathcal{F}}(\mathbb{Q}[SL_2(\mathbb{Z})]) \cong \mathbb{Z}[u_1, u_2, u_3] \Big/ \begin{array}{l} u_1^2 - 1, u_2^2 - u_2 - 2, \\ u_3^2 - 2u_1 - 2, \\ u_1 u_3 - u_3, \, u_1 u_2 - 2u_1 - u_2 + 2. \end{array}$$

The rank here is five, corresponding to the fact that $SL_2(\mathbb{Z})$ has precisely this number of conjugacy classes of cyclic subgroups, respectively isomorphic to $1$, $\mathbb{Z}/2$, $\mathbb{Z}/3$, $\mathbb{Z}/4$, and $\mathbb{Z}/6$.

EXAMPLE 3.2.

$$\Gamma = SL_3(\mathbb{Z}), \quad \mathcal{F} = \mathcal{C}(2), \text{ the family of all cyclic} \\ \text{subgroups of order a finite power of } 2.$$

Up to conjugacy, there are precisely two maximal subgroups in this family (both of order four), namely

$$S = \left\langle \begin{pmatrix} 1 & 0 & 1 \\ 0 & 0 & -1 \\ 0 & 1 & 0 \end{pmatrix} \right\rangle \quad \text{and} \quad T = \left\langle \begin{pmatrix} 0 & -1 & 0 \\ 1 & 0 & 0 \\ 0 & 0 & 1 \end{pmatrix} \right\rangle.$$

In addition, we have that $W_\Gamma(S) \cong W_\Gamma(T) \cong \mathbb{Z}/2$, and we have a short exact sequence

$$0 \longrightarrow G_{\mathcal{F}}(\mathbb{Q}[\Gamma]) \longrightarrow G(\mathbb{Q}[S])^{\mathbb{Z}/2} \oplus G(\mathbb{Q}[T])^{\mathbb{Z}/2} \longrightarrow \mathbb{Z} \longrightarrow 0,$$

as the two maximal classes intersect trivially. However, the action of $\mathbb{Z}/2$ is trivial on $G(\mathbb{Q}[\mathbb{Z}/4])$, and we obtain

$$G_{\mathcal{C}(2)}(\mathbb{Q}[SL_3(\mathbb{Z})]) \cong \mathbb{Z}[x_1, x_2, x_3, x_4] \Big/ \begin{array}{l} x_1^2 - 1, x_2^2 - 2x_1 - 2, \\ x_3^2 - 1, x_4^2 - 2x_3 - 2, \\ x_1 x_2 - x_2, x_1 x_3 - x_1 - x_3 + 1, \\ x_3 x_4 - x_4, x_2 x_4 - 2x_2 - 2x_4 + 4, \\ x_1 x_4 - 2x_1 - x_4 + 2, x_2 x_3 - 2x_3 - x_2 + 2. \end{array}$$

Note that the rank is five, corresponding to the five conjugacy classes $1, \langle S^2 \rangle, \langle T^2 \rangle, \langle S \rangle, \langle T \rangle$.

It is interesting to compare these calculations with the analogous ones over $\mathbb{C}$ (see [A]):

$$G_\mathcal{F}(\mathbb{C}[\mathbb{Z}/3 * \mathbb{Z}/5]) \cong \mathbb{Z}[z] \bigg/ z^7 + z^6 + z^5 - z^2 - z - 1,$$

$$G_\mathcal{F}(\mathbb{C}[SL_2(\mathbb{Z})]) \cong \mathbb{Z}[u] \bigg/ u^8 + u^6 - u^2 - 1,$$

$$G_{\mathcal{C}(2)}(\mathbb{C}[SL_3(\mathbb{Z})]) \cong G_{\mathcal{C}(2)}(\mathbb{Q}[SL_3(\mathbb{Z})]).$$

This last isomorphism is due to the fact that both elements of order 4 in $\langle S \rangle, \langle T \rangle$ are conjugate, and hence, conjugacy classes of cyclic 2–groups are in 1–1 correspondence with conjugacy classes of elements of order a finite power of two.

EXAMPLE 3.3. Let $\Gamma = GL_{p-1}(\mathbb{Z})$, $p$ an odd prime, and $\mathcal{F} = \mathcal{F}(p)$. From [CR] it is known that the only finite $p$–subgroups in $\Gamma$ are copies of $\mathbb{Z}/p$. In fact, if $\xi$ is a primitive $p$–th root of unity, and $R = \mathbb{Z}[\xi] \leq \mathbb{Q}(\xi)$ is the ring of algebraic integers, then these conjugacy classes of elements of order $p$ in $\Gamma$ correspond to the $R$–ideal classes in $\mathbb{Q}(\xi)$. Denote by I the set of ideal classes, of cardinality $Cl(p)$, the class number of $p$. Then the Galois group $\Delta \cong \mathbb{Z}/(p-1)$ acts on $I$, and it decomposes into orbits

$$I = \coprod_{i=1}^{t(p)} \Delta/S_i, \qquad S_i \subseteq \Delta.$$

Hence, we have a short exact sequence

$$0 \longrightarrow G_{\mathcal{F}(p)}(R[\Gamma]) \longrightarrow \left( \bigoplus_1^{Cl(p)} G(R[\mathbb{Z}/p]) \right)^\Delta \longrightarrow \mathbb{Z}^{t(p)-1} \longrightarrow 0$$

which becomes

$$0 \longrightarrow G_{\mathcal{F}(p)}(R[\Gamma]) \longrightarrow \bigoplus_1^{t(p)} G(R[\mathbb{Z}/p])^{S_i} \longrightarrow \mathbb{Z}^{t(p)-1} \longrightarrow 0.$$

Let $G_i = \mathbb{Z}/p \times_T S_i$ (the semi–direct product); the sequence above shows that

$$G_{\mathcal{F}(p)}(R[GL_{p-1}(\mathbb{Z})]) \cong G_{\mathcal{F}(p)}(R[G_1 * G_2 * \cdots * G_{t(p)}]).$$

(i) $R = \mathbb{C}$; in this case
$$G(\mathbb{C}[\mathbb{Z}/p]) \cong \mathbb{Z}[x] \bigg/ (x^p - 1)$$

and as a $\Delta$–module
$$G(\mathbb{C}[\mathbb{Z}/p]) \cong \mathbb{Z} \oplus \mathbb{Z}[\Delta]$$

hence,
$$G(\mathbb{C}[\mathbb{Z}/p])^{S_i} \cong \mathbb{Z}^{[\Delta:S_i]+1}$$

and $G_{\mathcal{F}(p)}(\mathbb{C}[GL_{p-1}(\mathbb{Z})])$ is free abelian of rank

$$\sum_{i=1}^{t(p)}(\mathbb{Z}^{[\Delta:S_i]}+1)-(t(p)-1)=Cl(p)+1$$

which is the expected result from Proposition 2.3.

(ii)  $R=\mathbb{Q}$; in this case

$$G(\mathbb{Q}[\mathbb{Z}/p])\cong \mathbb{Z}[w]\Big/w^2-(p-2)w-(p-1),$$

and $\Delta$ acts trivially on it.

Hence, the previous short exact sequence becomes

$$0\to G_{\mathcal{F}(p)}(\mathbb{Q}[GL_{p-1}(\mathbb{Z})])\longrightarrow \bigoplus^{t(p)}\mathbb{Z}[w]\Big/w^2-(p-2)w-(p-1)\longrightarrow \mathbb{Z}^{t(p)-1}\to 0,$$

and $G_{\mathcal{F}(p)}(\mathbb{Q}[GL_{p-1}(\mathbb{Z})])$ is free abelian of rank $2t(p)-t(p)+1=t(p)+1$ as was predicted by Proposition 2.2. Note that

$$G_{\mathcal{F}(p)}(\mathbb{Q}[GL_{p-1}(\mathbb{Z})])\cong G_{\mathcal{F}(p)}(\mathbb{Q}[\overset{t(p)}{*}\mathbb{Z}/p])$$

(iii)  $R=\mathbb{Z}$, in this case [Sw3]

$$G(\mathbb{Z}[\mathbb{Z}/p])\cong G(\mathbb{Q}[\mathbb{Z}/p])\times C_0(\mathbb{Z}[\mathbb{Z}/p])$$

where $C_0(\mathbb{Z}[\mathbb{Z}/p])$ is the reduced projective class group of $\mathbb{Z}[\mathbb{Z}/p]$, a finite group of order $Cl(p)$, and with product determined by the multiplication of ideals (this is in fact I). We obtain

$$\left|C_0(\mathbb{Z}[\mathbb{Z}/p])^{S_i}\right|=\prod_{j=1}^{t(p)}\frac{(p-1)|S_i\cap S_j|}{|S_i||S_j|}$$

and so there is an exact sequence

$$0\to A\to G_{\mathcal{F}(p)}(\mathbb{Z}[GL_{p-1}(\mathbb{Z})])\to G_{\mathcal{F}(p)}(\mathbb{Q}[\overset{t(p)}{*}\mathbb{Z}/p])\to 1$$

where $A$ is a finite group of order

$$|A|=\prod_{i=1}^{t(p)}\left[\prod_{j=1}^{t(p)}\frac{(p-1)|S_i\cap S_j|}{|S_i||S_j|}\right].$$

In particular, we obtain that $Cl(p)=1$, if and only if

$$G_{\mathcal{F}(p)}(\mathbb{Q}[GL_{p-1}(\mathbb{Z})])\cong G(\mathbb{Q}[\mathbb{Z}/p]).$$

EXAMPLE 3.4. We now consider a few of the so–called Bianchi groups. Let $K = \mathbb{Q}(\sqrt{d})$, with $d \in \mathbb{Z}$, $d < 0$ and squarefree, be an imaginary number field. We denote by $\mathcal{O}_d$ the ring of integers in $K$. As an abelian group $\mathcal{O}_d$ has as basis $\{1, w\}$, where

$$w = \begin{cases} \sqrt{d} & d \equiv 2, 3 \mod 4 \\ \frac{1+\sqrt{d}}{2} & d \equiv 1 \mod 4. \end{cases}$$

The Bianchi group $\Gamma_d$ is defined as

$$\Gamma_d = PSL_2(\mathcal{O}_d).$$

These groups satisfy the FCP–condition. For small values of $d$, Flöge [F] has determined the lattice of finite subgroups (see also [G-S]). The elements of finite order in these groups are of order 1, 2 or 3. Hence we use the families $\mathcal{F}(2)$ and $\mathcal{F}(3)$.

(1) $\Gamma = \Gamma_{-1} = PSL_2(\mathbb{Z}[\sqrt{-1}])$.

The finite subgroups (up to conjugacy) are determined by the diagram of intersections:

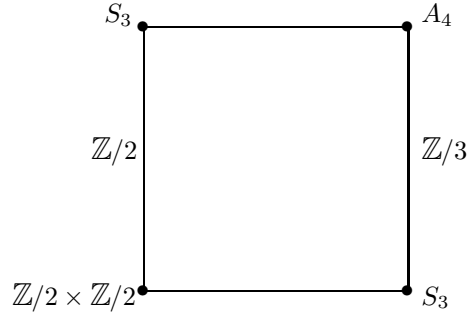

In fact, this simply describes an isomorphism

$$\Gamma_{-1} \cong ((\mathbb{Z}/2)^2 \underset{\mathbb{Z}/2}{*} S_3) \underset{\mathbb{Z}/2 * \mathbb{Z}/3}{*} (A_4 \underset{\mathbb{Z}/3}{*} S_3).$$

Let $\Gamma_I = (\mathbb{Z}/2)^2 \underset{\mathbb{Z}/2}{*} S_3$, $\Gamma_{II} = A_4 \underset{\mathbb{Z}/3}{*} S_3$; recall that

$$G_{\mathcal{F}(2)}(\mathbb{C}[(\mathbb{Z}/2)^2]) \cong \mathbb{Z}[x, y] \Big/ \begin{array}{l} x^2 - 1, \\ y^2 - 1, \end{array}$$

$$G_{\mathcal{F}(2)}(\mathbb{C}[S_3]) \cong \mathbb{Z}[z] \Big/ z^2 - 1 \cong G_{\mathcal{F}(2)}(\mathbb{C}[\mathbb{Z}/2]),$$

from which we deduce

$$G_{\mathcal{F}(2)}(\mathbb{C}[\Gamma_I]) \cong \mathbb{Z}[x,y] \bigg/ \begin{array}{l} x^2-1, \\ y^2-1, \end{array}$$

and

$$G_{\mathcal{F}(3)}(\mathbb{C}[\Gamma_I]) \cong G_{\mathcal{F}(3)}(\mathbb{C}[S_3]) \cong \mathbb{Z}[\gamma] \bigg/ \gamma^2-\gamma-2\,.$$

Similarly,

$$G_{\mathcal{F}(3)}(\mathbb{C}[A_4]) \cong \mathbb{Z}[w] \bigg/ w^3-1\,,$$

and

$$G_{\mathcal{F}(2)}(\mathbb{C}[A_4]) \cong \left( \mathbb{Z}[x,y] \bigg/ \begin{array}{l} x^2-1, \\ y^2-1 \end{array} \right)^{\mathbb{Z}/3},$$

where $\mathbb{Z}/3$ acts by $x \longmapsto y \longmapsto xy$, fixing $1$.
Hence,

$$G_{\mathcal{F}(2)}(\mathbb{C}[A_4]) \cong \mathbb{Z}[y] \bigg/ y^2-2y-3\,.$$

We deduce that

$$G_{\mathcal{F}(3)}(\mathbb{C}[\Gamma_{II}]) \cong G_{\mathcal{F}(3)}(\mathbb{C}[S_3]) \cong \mathbb{Z}[\gamma] \bigg/ \gamma^2-\gamma-2\,,$$

and

$$G_{\mathcal{F}(2)}(\mathbb{C}[\Gamma_{II}]) \cong \mathbb{Z}[u_1, u_2] \bigg/ \begin{array}{l} u_1^2-1, \\ u_2^2-2u_2-3\,. \end{array}$$

Similarly, we have

$$G_{\mathcal{F}(2)}(\mathbb{C}[\mathbb{Z}/2 * \mathbb{Z}/3]) \cong G(\mathbb{C}[\mathbb{Z}/2])\,, \qquad G_{\mathcal{F}(3)}(\mathbb{C}[\mathbb{Z}/2 * \mathbb{Z}/3]) \cong G(\mathbb{C}[\mathbb{Z}/3])\,.$$

At $p=2$ we obtain an exact sequence:

$$0 \to G_{\mathcal{F}(2)}(\mathbb{C}[\Gamma_{-1}]) \longrightarrow \mathbb{Z}[x,y] \bigg/ \begin{array}{l} x^2-1, \\ y^2-1 \end{array} \bigoplus \mathbb{Z}[u_1, u_2] \bigg/ \begin{array}{l} u_1^2-1, \\ u_2^2-2u_2-3 \end{array} \xrightarrow{\phi} \mathbb{Z}[z] \bigg/ z^2-1 \to 0\,,$$

where

$$\phi(x) = z\,, \qquad \phi(y) = 1\,, \qquad \phi(u_1) = -z\,, \qquad \phi(u_2) = -3\,.$$

¿From this, we obtain

$$\boxed{G_{\mathcal{F}(2)}(\mathbb{C}[\Gamma_{-1}]) \cong \mathbb{Z}[v_1, v_2, v_3] \bigg/ \begin{array}{l} v_1^2-1, v_2^2-1 \\ v_3^2-2v_3-3, \\ v_2 v_3 - 3v_2 - v_3 - 3 \end{array}}$$

At $p=3$ we have an exact sequence

$$0 \to G_{\mathcal{F}(3)}(\mathcal{C}[\Gamma_{-1}]) \longrightarrow \mathbb{Z}[\gamma_1]\Big/\gamma_1^2-\gamma_1-2 \bigoplus \mathbb{Z}[\gamma_2]\Big/\gamma_2^2-\gamma_2-2 \xrightarrow{\Psi} \left(\mathbb{Z}[w]\Big/w^3-1\right)^{\mathbb{Z}/2}.$$

In this case, note that $\psi$ maps to the $\mathbb{Z}/2$ invariants, as $W_\Gamma(\mathbb{Z}/3)=\mathbb{Z}/2$. Of course, these invariants are isomorphic to $G_{\mathcal{F}(3)}(\mathbb{C}[S_3])$, and we obtain

$$\boxed{G_{\mathcal{F}(3)}(\mathbb{C}[\Gamma_{-1}]) \cong \mathbb{Z}[\gamma]\Big/\gamma^2-\gamma-2\,.}$$

We conclude that $n_2(\Gamma_{-1})=5$, $n_3(\Gamma_{-1})=2$, and $n(\Gamma_{-1})=6$.

(2) $\Gamma_d = \Gamma_{-3} = PSL_2(\mathbb{Z}[\frac{1+\sqrt{-3}}{2}])$.

In this case, the finite subgroup structure is described by the diagram

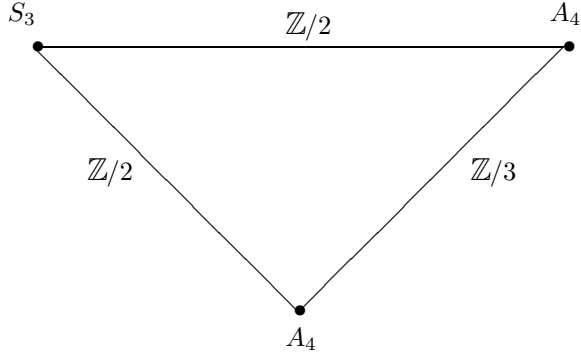

At $p=2$ this yields a simplified exact sequence

$$0 \longrightarrow G_{\mathcal{F}(2)}(\mathcal{C}[\Gamma_{-3}]) \longrightarrow \mathbb{Z}[y_1]\Big/y_1^2-2y_1-3 \bigoplus \mathbb{Z}[y_2]\Big/y_2^2-2y_2-3 \xrightarrow{\rho} \mathbb{Z}[z]\Big/z^2-1 \longrightarrow 0,$$

where $\rho(y_1)=2z+1$, $\rho(y_2)=-2z-1$.
¿From this we obtain

$$\boxed{G_{\mathcal{F}(2)}(\mathbb{C}[\Gamma_{-3}]) \cong \mathbb{Z}[y]\Big/y^2-2y-3}$$

Similarly, for $p=3$, we get

$$0 \longrightarrow G_{\mathcal{F}(3)}(\mathcal{C}\Gamma_{-3}) \longrightarrow \mathbb{Z}[w]\Big/w^3-1 \bigoplus \mathbb{Z}[\gamma]\Big/\gamma^2-\gamma-2 \longrightarrow \mathbb{Z} \longrightarrow 0,$$

from which

$$\boxed{G_{\mathcal{F}(3)}(\mathbb{C}[\Gamma_{-3}]) \cong \mathbb{Z}[\alpha,\beta]\Big/\begin{array}{l}\alpha^3-1,\beta^2-\beta-2,\\ \alpha^2\beta-2\alpha^2-\beta+2\end{array}}$$

We conclude that $n_2(\Gamma_{-3}) = 2$, $n_3(\Gamma_{-3}) = 4$, and $n(\Gamma) = 5$.

These examples illustrate the *computability* of $G_{\mathcal{F}}(\mathbb{K}[\Gamma])$ for $\mathbb{K} = \mathbb{Q}$ or $\mathbb{C}$.

## 4. Projective Modules

In this section we will relate our construction to $K_0(\mathbb{K}[\Gamma])$, the Grothendieck group of projective modules, $\mathbb{K}$ a field of characteristic $0$. Denote by $K_\Gamma^0(\mathbb{K})$ the Grothendieck group corresponding to the category of $\mathbb{K}[\Gamma]$– modules which are finitely generated over $\mathbb{K}$ with relations induced only by direct sum decompositions. Clearly, the group $G^f(\mathbb{K}[\Gamma])$ is a quotient of this.

We recall that Loday [L] has constructed a bilinear pairing, natural in $\Gamma$

$$(4.1) \qquad K_\Gamma^0(\mathbb{K}) \times K_0(\mathbb{K}[\Gamma]) \xrightarrow{\varepsilon_\Gamma} K_0(\mathbb{K}).$$

We briefly describe how it is defined: if $[V] \in K_\Gamma^0(\mathbb{K})$, where $V$ is determined by a homomorphism $\Gamma \xrightarrow{\rho} GL_n(\mathbb{K})$, then it extends to a map of algebras

$$\rho^{ann} : \mathbb{K}[\Gamma] \longrightarrow M_{n \times n}(\mathbb{K}),$$

and hence, induces

$$\rho_*^{ann} : K_0(\mathbb{K}[\Gamma]) \longrightarrow K_0(M_{n \times n}(\mathbb{K})) \cong K_0(\mathbb{K}).$$

Then,

$$\varepsilon([V], x) = \rho_*^{ann}(x);$$

this defines $\varepsilon$ on generators, and it extends by linearity.

If $Q \subseteq \Gamma$ is any subgroup, then there is a well– defined induction map

$$K_0(\mathbb{K}[Q]) \xrightarrow{\mathrm{Ind}_Q^\Gamma} K_0(\mathbb{K}[\Gamma])$$

determined on the generators by

$$M \longmapsto M \bigotimes_Q \mathbb{K}[\Gamma].$$

Apply the functor $\mathrm{Hom}(-, \mathbb{Z})$ to this, we get:

$$\mathrm{Hom}(K_0(\mathbb{K}[\Gamma]), \mathbb{Z}) \xrightarrow{(\mathrm{Ind}_Q^\Gamma)^*} \mathrm{Hom}(K_0(\mathbb{K}[Q]), \mathbb{Z}).$$

Now, consider $\mathcal{F}$, a family of finite subgroups, and suppose $\Gamma$ is FCP. Then, as before, we can define

$$(4.2) \qquad \mathrm{Hom}(K_0(\mathbb{K}[\Gamma]), \mathbb{Z})_{\mathcal{F}} = \varprojlim_{H \in \mathcal{F}} \mathrm{Hom}(K_0(\mathbb{K}[H]), \mathbb{Z}),$$

using conjugation and the $\mathrm{Ind}^*$.

The following diagram is commutative:

$$
\begin{array}{ccc}
K^0_\Gamma(\mathbb{K}) & \longrightarrow & \mathrm{Hom}(K_0(\mathbb{K}[\Gamma]),\,\mathbb{Z}) \\
\downarrow \mathrm{res} & & \downarrow \mathrm{Ind}^* \\
K^0_H(\mathbb{K}) & \longrightarrow & \mathrm{Hom}(K_0(\mathbb{K}[H]),\,\mathbb{Z})
\end{array}
$$
(4.3)

which is the extension of the usual adjunction formula for finite groups. In particular, the map on the bottom is induced by the usual pairing on character rings, which coincides with $\varepsilon_H$. Tensoring with $\mathbb{K}$, and taking limits, we obtain the commutative diagram

$$
\begin{array}{ccc}
K^0_\Gamma(\mathbb{K}) \otimes \mathbb{K} & \longrightarrow & \mathrm{Hom}(K_0(\mathbb{K}[\Gamma]),\,\mathbb{K}) \\
\downarrow \mathrm{res} & & \downarrow \mathrm{Ind}^*_\mathcal{F} \\
G_\mathcal{F}(\mathbb{K}[\Gamma]) \otimes \mathbb{K} & \longrightarrow & \mathrm{Hom}(K_0(\mathbb{K}[\Gamma]),\mathbb{K})_\mathcal{F}
\end{array}
$$
(4.4)

Note that $\mathrm{I}nd$ is the formal adjoint of the restriction for finite groups, hence, the bottom map is an isomorphism. We obtain (using the notation in section 2)

PROPOSITION 4.1. *Under the above hypotheses, if $\Gamma$ is VTF then*

$$dim_\mathbb{Q} K_0(\mathbb{K}[\Gamma]) \otimes \mathbb{Q} \geq dim_\mathbb{Q}\,(imj_\mathcal{F}) \otimes \mathbb{Q} \qquad \square$$

In particular if $\Gamma$ is VTF and $\mathbb{K} = \mathbb{C}$, then by Proposition 2.19 the rank of this image is bounded below by $\mathrm{Max}(n_G)$, as $G$ ranges over all finite quotients of $\Gamma$, and as a corollary of this we obtain: if $K_0(\mathbb{C}[\Gamma]) \cong \mathbb{Z}$, then $\Gamma$ is torsion free. This is the converse of a well–known conjecture in algebraic K–theory (see [F–J]). It would be interesting to determine when the vertical arrow on the right hand side of (4.4) is an isomorphism.

## 5. Equivariant K–theory

Families of finite subgroups play a natural role in representation theory. ¿From this it is logical to expect that they also will be important in the topological invariants arising from vector bundles, in particular for K–theory. In this section we will apply families of finite subgroups to prove a general result about the equivariant K–theory of a finite $G$–CW complex, where $G$ is a finite group. However this will have an interesting application to discrete groups of the type we have been considering. This application was implicitly used in [A], and a similar analysis has been done by Quillen [Q] for equivariant cohomology.

In what follows we will assume that $Y$ is a finite $G$–CW complex, and that $G$ is a finite group. Consider the collection of pairs $(H, c)$ where $c \in \pi_0(Y^H)$, $H \subset G$ a subgroup. We define a morphism $\theta : (H, c) \to (H', c')$ as a triple $((H, c), (H', c'), \overline{\theta})$ where $\overline{\theta}$ is a homomorphism from $H$ to $H'$ of the form $\overline{\theta}(h) = ghg^{-1}$, for some $g \in G$ such that $gHg^{-1} \subset H'$, $c' \subset gc$. In this way we may construct a category (see [Q]) and in particular form the limit of complex representation rings

$$\mathcal{R}_G(Y) = \lim_{(H,c)} R(H).$$

The map $(H, pt) \to (G, Y)$ consisting of the inclusion homomorphism from $H$ to $G$ and the map from the one–point space to $Y$, taking the point to the component $c \in Y^H$ induces a map of rings $K_G^0(Y) \to R(H)$ and hence a ring homomorphism to the compatible elements,

$$\varphi : K_G^0(Y) \to \mathcal{R}_G(Y).$$

We have

PROPOSITION 5.1. *Let $Y$ be a finite $G$–CW complex, then the map $\varphi : K_G^0(Y) \longrightarrow \mathcal{R}_G(Y)$ has nilpotent kernel and becomes a surjection after tensoring with the rationals.*

PROOF. The proof involves adapting Quillen's methods to equivariant K–theory. We describe how this goes. There is a spectral sequence in equivariant K–theory [Seg], with

$$E_2^{p,q} \cong H^p(Y/G, \mathcal{K}_G^q) \to K_G^{p+q}(Y),$$

where $\mathcal{K}_G^t$ is the sheaf on $Y/G$ associated to the presheaf $V \mapsto K_G^t(q^{-1}V)$, and $q : Y \to Y/G$ is the natural projection.

Let $s \in H^0(Y/G, \mathcal{K}_G^0)$, and $K \subset G$ a subgroup. Let $s_K : Y^K \to R(K)$ be the function whose value at $y$ is the image of $s(qy) \in K_G^0(Gy)$ under the homomorphism $K_G^0(Gy) \to K_K^0(pt)$ associated to the inclusion of $K$ in $G$ and the map $pt \to Gy$ with image $\{y\}$. This function is locally constant, and the family $\{s_K\}$ is compatible. Hence we have defined a homomorphism

$$H^0(Y/G, \mathcal{K}_G^0) \to \mathcal{R}_G(Y)$$

which is in fact an isomorphism, as can be verified by looking at the stalks of this map of sheaves (see [Q], Prop.5.11).

The edge homomorphism in the spectral sequence gives rise to a ring homomorphism

$$K_G^0(Y) \to H^0(Y/G, \mathcal{K}_G^0)$$

which has nilpotent kernel ([Q] Prop.3.2). This can be identified with the map $\varphi$, proving the first part of the proposition.

To prove the rest, we first recall how the following additive decomposition was proved in [K]:

$$K_G^*(Y) \otimes \mathbb{C} \cong \bigoplus_{(g)} K^*(Y^g)^{C(g)} \otimes \mathbb{C}$$

where $(g)$ runs over conjugacy classes of elements in $G$, $Y^g$ is the fixed–point set of the action of $g \in G$, and $C(g)$ is the centralizer of $g$ in $G$. The explicit isomorphism is described a follows: denote $H =<g>$, then $K_H^*(Y^H) \cong R(H) \otimes K^*(Y^H)$; we have a restriction map $K_G^*(Y) \to K_H^*(Y^H)$ as well as a map $K_H^*(Y^H) \to K^*(Y^H) \otimes \mathbb{C}$ defined by mapping the generator of $R(H)$ to the appropriate finite root of unity. The sum of these compositions lands in the term on the right hand side of the formula above, and becomes an isomorphism after complexification. At the level of path components we obtain a surjection $K_G^0(Y) \otimes \mathbb{C} \to (\mathbb{C})^M$, where $M$ is precisely the sum over all conjugacy classes of $g \in G$ of the number of distinct $C(g)$–orbits in $\pi_o(Y^g)$. From the description above, it is clear that this surjection factors through the complexification of $\varphi$.

On the other hand consider the $G$–action induced by conjugation on the set

$$S_G(Y) = \{(g,c) \mid c \in \pi_0(Y^g)\},$$

defined by $(g,c) \mapsto (hgh^{-1}, hc)$. Using characters we can identify $\mathcal{R}_G(Y) \otimes \mathbb{C}$ with the ring of complex valued class functions on $S_G(Y)$. From this it follows that the rank of $\mathcal{R}_G(Y)$ is equal to $M$. We infer that $\varphi$ induces a rational surjection. □

Now we apply this in the context of discrete groups. First we introduce

DEFINITION 5.1. *If $\Gamma$ is a discrete group, and $\mathcal{F}$ is a family of finite subgroups, then an $\mathcal{F}$–admissible $\Gamma$–complex is a finite dimensional $\Gamma$–CW complex $X$ satisfying the following conditions:*
 *(i) $X^H \neq \emptyset$ if and only if $H \in \mathcal{F}$*
 *(ii) $X^H$ is contractible for all $H \in \mathcal{F}$*
 *(iii) $X/\Gamma$ is compact.*

Note: if $\mathcal{F}$ is the family of all finite subgroups in $\Gamma$, $X$ is simply said to be admissible.

Let us assume for the remainder of this section that $\Gamma$ is a VTF group and that for a fixed family $\mathcal{F}$ of finite subgroups there exists an $\mathcal{F}$–admissible $\Gamma$–complex $X$. Choose a torsion–free normal subgroup $\Gamma' \subset \Gamma$, with quotient a finite group $G$. An obvious consequence of our assumptions is that $\Gamma'$ acts freely on $X$. A less apparent consequence is that $\mathcal{F}$ will have a finite number of elements up to conjugacy; in particular if $\mathcal{F}$ is the family of all finite subgroups in $\Gamma$, then the group will satisfy the FCP condition. For these and related facts we refer to [Q] and [A].

REMARK Admissible complexes exist for many interesting classes of groups, including arithmetic and mapping class groups. It is also possible to work $p$–locally for these groups (see [JMO]).

We now apply the preceding proposition to $Y = X/\Gamma'$, $G = \Gamma/\Gamma'$. Using that $X^H \neq \emptyset$ all $H \in \mathcal{F}$, the restriction maps followed by fixed–point inclusions

$$K_\Gamma^*(X) \longrightarrow K_H^*(X) \longrightarrow G(\mathbb{C}[H])$$

for each $H \in \mathcal{F}$, induce a map

$$K_\Gamma^0(X) \xrightarrow{\varphi_\mathcal{F}} G_\mathcal{F}(\mathbb{C}[\Gamma]),$$

as it is compatible with respect to restriction and conjugation. Our hypotheses imply that $\mathcal{R}_\Gamma(X)$ is well-defined and naturally isomorphic to $G_\mathcal{F}(\mathbb{C}[\Gamma])$. In fact the map $(X, \Gamma) \to (Y, G)$ gives rise to a commutative diagram

$$\begin{array}{ccc} K_\Gamma^*(X) & \xrightarrow{\varphi_\mathcal{F}} & G_\mathcal{F}(\mathbb{C}[\Gamma]) \\ \uparrow & & \uparrow \\ K_G^*(Y) & \xrightarrow{\varphi} & \mathcal{R}_G(Y) \end{array}.$$

As pointed out by Quillen [Q], the vertical arrow on the right is an isomorphism; as $\Gamma'$ acts freely on $X$, the one on the left is also one, and hence we obtain

COROLLARY 5.1. *If $X$ is an $\mathcal{F}$-admissible $\Gamma$-complex, then there is a ring map $K_\Gamma^0(X) \xrightarrow{\varphi_\mathcal{F}} G_\mathcal{F}(\mathbb{C}[\Gamma])$ which has nilpotent kernel and is a rational surjection.*

This result shows that the equivariant K–theory of $X$ can be used to distinguish conjugacy classes of elements of finite order in $\Gamma$, *unlike the actual representation ring*, as we mentioned before. Likewise we see that our construction is a basic ingredient in equivariant K–theory. We remark that it would be both natural and interesting to introduce Adams operations in this context.

Next we provide an application of our methods to certain virtually free groups. Let $\Gamma$ be a group (of virtual cohomological dimension 1) which has a tree $T$ as an admissible complex. By [Se1] there exist finite subgroups $G_1, \ldots, G_m$ and $H_1, \ldots, H_n$ such that $T/\Gamma$ has precisely m orbit representatives of 0-simplices, with respective isotropy $G_i$ and n orbit representatives of 1-simplices, with respective isotropy $H_j$. The $G_i$'s are maximal finite subgroups, and the $H_j$'s correspond to their intersections. Note that every finite subgroup can be conjugated into a subgroup of these maximal ones. In this case the spectral sequence for equivariant K–theory simplifies to yield:

$$0 \longrightarrow K_\Gamma^0(T) \longrightarrow \oplus_{i=1}^m R(G_i) \longrightarrow \oplus_{j=1}^n R(H_j) \longrightarrow K_\Gamma^1(T) \longrightarrow 0$$

from which we conclude

$$K_\Gamma^0(T) \cong G_\mathcal{F}(\mathbb{C}[\Gamma]).$$

Rationally it turns out [A] that in fact $T$ can be chosen so that

$$K_\Gamma^1(T) \otimes \mathbb{Q} \cong \oplus_{(\gamma)} H^1(BC_\Gamma(\gamma), \mathbb{Q})$$

where the sum is taken over conjugacy classes of elements of finite order in $\Gamma$. We obtain, using the notation from section 1:

COROLLARY 5.2. *The number of conjugacy classes of elements of finite order in a group $\Gamma$ having a tree as an admissible complex can be computed from the formula*

$$n(\Gamma) = n(G_1) + \cdots + n(G_m) - (n(H_1) + \cdots + n(H_n)) + \sum_{(\gamma)} dim_{\mathbb{Q}} H^1(BC_\Gamma(\gamma), \mathbb{Q}).$$

$\square$

We conclude by recording a few questions related to the topics in this paper (the relevant section is mentioned after each):

1. If $\mathbb{K}$ is the field of fractions of $R$, a Noetherian integral domain, calculate the map $G_\mathcal{F}(R[\Gamma]) \longrightarrow G_\mathcal{F}(\mathbb{K}[\Gamma])$. [section 2]

2. Under what hypotheses on $\Gamma$ is the inclusion $\overline{G}_\mathcal{F}(R[\Gamma]) \longrightarrow G_\mathcal{F}(R[\Gamma])$ a rational equivalence? [section 2]

3. For what groups is is true that

$$Hom(K_0(\mathbb{K}[\Gamma]), \mathbb{K}) \longrightarrow Hom(K_0(\mathbb{K}[\Gamma]), \mathbb{K})_\mathcal{F}$$

is an isomorphism? [section 4]

4. If $\mathbb{F}$ is a suitable field of characteristic $p > 0$, find a formula for $K\mathbb{F}^0(B\Gamma)$ (modular K–theory) using $G_{\mathcal{F}(p')}(\mathbb{C}[\Gamma])$ and Adams operations. [section 5]

Department of Mathematics, University of Wisconsin, Madison, Wisconsin, 53706
*E-mail address*: adem@math.wisc.edu